\documentclass[12pt,a4paper]{article}

\usepackage[latin2]{inputenc}
\usepackage{amssymb}
\usepackage{paralist}
\usepackage{amsmath}
\usepackage{amsfonts}
\usepackage{amsthm}
\usepackage{fullpage}
\usepackage{enumerate}
\newtheorem{thm}{Theorem}

\newtheorem{prob}{Problem}
\newtheorem{conj}{Conjecture}
\newtheorem{cor}{Corollary}
\newtheorem{lem}{Lemma}

\begin{document}

\title{Partitions of the set of nonnegative integers with the same
representation functions}
\author{S\'andor Z. Kiss \thanks{Institute of Mathematics, Budapest
University of Technology and Economics, H-1529 B.O. Box, Hungary,
kisspest@cs.elte.hu
This author was supported by the OTKA Grant No. NK105645. This research was 
supported by the
National Research, Development and Innovation Office NKFIH Grant No. K115288.}, 
Csaba S\'andor \thanks{Institute of Mathematics, Budapest University of
Technology and Economics, H-1529 B.O. Box, Hungary, csandor@math.bme.hu.
This author was supported by the OTKA Grant No. K109789. This paper was
supported by the J\'anos Bolyai Research Scholarship of the Hungarian
Academy of Sciences.}
}
\date{}
\maketitle

\begin{abstract}
For a set of nonnegative integers $S$ let $R_{S}(n)$ denote the number
 of unordered representations of the integer $n$ as the sum of two
 different terms from $S$. In this paper we focus on partitions of the
  natural numbers into two sets affording identical 
 representation functions. We solve a recent problem of Lev
 and Chen.
\end{abstract}
\textit{2000 AMS \ Mathematics subject classification number}: 11B34.

\section{Introduction}

Let $S$ be a set of nonnegative integers and let $R_{S}(n)$ denote the number
 of solutions of the equation $s + s^{'} = n$, where $s$, $s^{'} \in S$ and $s
 < s^{'}$. Let $\mathbb{N}$ denotes the set of nonnegative integers.
The binary representation of an integer $n$ is the representation of
$n$ in the number system with base $2$. Let $A$ be the set of those
nonnegative integers which contains even number of $1$ binary digits in its
binary representation and let $B$ be the complement of $A$. The set $A$ is 
called Thue-Morse sequence.
The investigation of the partitions of the set of nonnegative
 integers with identical representation functions was a popular topic in the
 last few decades [1], [2], [7], [8], [9]. By using the Thue - Morse sequence 
in 2002  Dombi [5] constructed two sets of nonnegative integers with infinite 
symmetric difference such that the corresponding representation functions are
 identical. Namely, he proved the following theorem.

\begin{thm}(Dombi, 2002)
The set of positive integers can be partitioned into two subsets $C$ and $D$ 
such that $R_{C}(n) = R_{D}(n)$ for every positive integer $n$.
\end{thm}

\noindent The complete description of the suitable partitions is the following.

\begin{thm}
Let $C$ and $D$ be sets of nonnegative integers such that 
$C\cup D=\mathbb{N}$, $C\cap D=\emptyset $ and $0\in C$. Then $R_C(n)=R_D(n)$ 
for every positive integer $n$ if and only if $C=A$ and $D=B$. 
\end{thm}

As far as we know this theorem has never been formulated in the above form, 
but the nontrivial part was proved by Dombi, therefore this theorem is only a 
little extension of Dombi's result. We give an alternative proof of the 
previous theorem.

A finite version of the above theorem is the following.
Put $A_{l} = A
\cap [0, 2^{l} - 1]$ and $B_{l} = B \cap [0, 2^{l} - 1]$.

\begin{thm}
Let $C$ and $D$ be sets of
nonnegative integers such that $C \cup D = [0, m]$ and
$C \cap D = \emptyset$, $0 \in C$. Then $R_{C}(n) = R_{D}(n)$ for every 
positive integer $n$ if and only if
there exists an $l$ natural number such that $C = A_{l}$ and $D = B_{l}$.
\end{thm}

\noindent If $C=A\cap [0,m]$ and $D=B\cap [0,m]$ then by Theorem 2 we have 
$R_{C}(n)=R_{D}(n)$ for $n\leq m$, therefore Theorem 3. implies the following 
corollary.

\begin{cor}
If $C = A \cap [0, m]$ and $D = B \cap [0, m]$, where $m$ is a positive 
integer  not of the form $2^l-1$, then there exists a positive integer 
$m < n < 2m$  such that $R_{C}(n) \ne R_{D}(n)$. 
\end{cor}

\noindent In Dombi's example the union of the set $C$ and $D$ is the set of nonnegative integers, and they are disjoint sets. Tang and Yu [11] proved that if the union of the sets $C$ and $D$ is the set of nonnegative integers and the representation functions are identical from a certain point on, then at least one cannot have the intersection of the two sets is the non-negative integers divisible by $4$ i.e.,
\begin{thm}(Tang and Yu, 2012)
If $C \cup D = \mathbb{N}$ and $C \cap D = 4\mathbb{N}$, then
$R_{C}(n) \ne R_{D}(n)$ for infinitely many $n$.
\end{thm}
\noindent Moreover, they conjectured that under the same
assumptions the intersection cannot be the union of infinite arithmetic 
progressions.

\begin{conj}(Tang and Yu, 2012)
Let $m \in \mathbb{N}$ and $R \subset \{0, 1, \dots{}, m - 1\}$. If $C \cup D = \mathbb{N}$ and $C \cap D = \{r + km: k \in \mathbb{N}, r \in R\}$, then
$R_{C}(n) = R_{D}(n)$ cannot hold for all sufficiently large $n$.
\end{conj}

\noindent Recently Tang extended Theorem 4. In particular, she proved the following theorem [10]:
\begin{thm}(Tang, 2016)
Let $k \ge 2$ be an integer. If $C \cup D = \mathbb{N}$ and $C \cap D = k\mathbb{N}$, then
$R_{C}(n) \ne R_{D}(n)$ for infinitely many $n$.
\end{thm}

\noindent Chen and Lev [2] disproved the above conjecture
by constructing a family of partitions of the set of natural numbers such that
all the corresponding representation functions are the same and the
intersection of the two sets is an infinite arithmetic progression properly
contained in the set of natural numbers.

\begin{thm}(Chen and Lev, 2016)
Let $l$ be a positive integer. There exist sets $C$ and $D$ such that 
$C \cup D = \mathbb{N}$, $C \cap D =
(2^{2l} - 1) + (2^{2l+1} - 1)\mathbb{N}$ and $R_{C}(n) = R_{D}(n)$ for every 
positive integer $n$.
\end{thm}

\noindent Their construction is based on the following lemma:

\begin{lem}(Chen and Lev, 2016)
If there exist sets $C_0$ and $D_0$ such that $C_0\cup D_0=[0,m]$, 
$C_0\cap D_0=\{ r\}$ and $R_{C_0}(n)=R_{D_0}(n)$ for every positive integer $n$ 
then there 
exist sets $C$ and $D$ such that $C\cup D=\mathbb{N}$, $C\cap D=r+(m+1)\mathbb{N}$ and $R_C(n)=R_D(n)$ for every positive integer $n$.
\end{lem}

\noindent When $r = 2^{2l} - 1$ and $m = 2^{2l+1} - 2$ they proved the existence of two such sets $C_{0}$ and $D_{0}$, which verifies the sufficiency of the following problem in [2] (we use different notations then they applied).

\begin{prob}(Chen and Lev, 2016)
Given $R_C(n)=R_D(n)$ for every positive integer $n$, $C\cup D=[0,m]$, and $C\cap D=\{ r\}$ with integers $r\geq 0$ and $m\geq 2$, must there exist an integer $l\geq 1$ such that $r=2^{2l}-1$, $m=2^{2l+1}-2$, $C = A_{2l} \cup (2^{2l} - 1 +
B_{2l})$ and $D = B_{2l} \cup (2^{2l} - 1 + A_{2l})$?
\end{prob}

\noindent In this paper we solve this problem affirmatively.

\begin{thm}
Let $C$ and $D$ be sets of
nonnegative integers such that $C \cup D = [0, m]$ and
$C \cap D = \{r\}$, $0 \in C$. Then $R_{C}(n) = R_{D}(n)$ for every $n$ positive 
integer if and only if
there exists an $l$ natural number such that $C = A_{2l} \cup (2^{2l} - 1 +
B_{2l})$ and $D = B_{2l} \cup (2^{2l} - 1 + A_{2l})$.
\end{thm}

\noindent The previous theorem suggests that there are no other counterexample for Tang and Yu's conjecture.

\begin{prob}(Chen and Lev, 2016)
Given $R_C(n)=R_D(n)$ for every $n$ positive integer, $C\cup D=\mathbb{N}$, and $C\cap D=r+m\mathbb{N}$ with integers $r\geq 0$ and $m\geq 2$, must there exist an integer $l\geq 1$ such that $r=2^{2l}-1$, $m=2^{2l+1}-1$?
\end{prob}

%
Similar questions were investigated for unordered representation functions in 
[3], [4], [6]. 
Thoughout this paper the characteristic function of the set $S$ is denoted by 
$\chi_{S}(n)$,
i.e.,
\[
\chi _S(n) = \left\{
\begin{aligned}
1 \textnormal{, if } n \in S \\
0 \textnormal{, if } n \notin S
\end{aligned} \hspace*{3mm}.
\right.
\]
For a nonnegative integer $a$ and a set of nonnegative integers $S$ we define
the sumset $a + S$ by
\[
a + S = \{a+b: b \in S\}.
\]

\section{Proof of Theorem 2. and 3.}

First we prove that if there exists a natural number $l$ such that $C =
A_{l}$ and $D = B_{l}$, then $R_{C}(n) = R_{D}(n)$ for every positive integer $n$.

We prove by induction on $l$. For $l = 1$, $A_{1} = \{0\}$ and $B_{1} =
\{1\}$ thus $R_{A_{1}}(n) = R_{B_{1}}(n)=0$. Assume the statement holds
for any $l$ and we prove it to $l + 1$. By the definition of $A$ and $B$ 
we have $A_{l+1}=A_l\cup (2^l+B_l)$ and $B_{l+1}=B_l\cup (2^l+A_l)$. Hence

\[
R_{A_{l+1}}(n) = R_{A_{l} \cup (2^{l} + B_{l})}(n) = |\{(a,a'): a <
a', a, a' \in A_{l}, a + a' = n\}|
\]
\[
+ |\{(a,a'): a \in A_{l},
a' \in 2^{l} + B_{l}, a + a' =
n\}| + |\{(a,a'): a < a', a, a' \in 2^{l} + B_{l}, a + a' =
n\}|
\]
\[
= R_{A_{l}}(n) + |\{(a,a'): a \in A_{l},
a' \in B_{l}, a + a' = n - 2^{l}\}| + R_{B_{l}}(n - 2^{l+1}).
\]
On the other hand

\[
R_{B_{l+1}}(n) = R_{B_{l} \cup (2^{l} + A_{l})}(n) = |\{(a,a'): a <
a', a, a' \in B_{l}, a + a' = n\}|
\]
\[
+ |\{(a,a'): a \in B_{l},
a' \in 2^{l} + A_{l}, a + a' =
n\}| + |\{(a,a'): a<a' a, a' \in 2^{l} + A_{l}, a + a' =
n\}|
\]
\[
= R_{B_{l}}(n) + |\{(a,a'): a \in B_{l},
a' \in A_{l}, a + a' = n - 2^{l}\}| + R_{A_{l}}(n - 2^{l+1}),
\]
thus we get the result.

Observe that if $k
\le 2^{l} - 1$, then $R_{A_{l}}(k) =  R_{A}(k)$ and $R_{B_{l}}(k) =
R_{B}(k)$. On the other hand $R_{A_{l}}(k) = R_{B_{l}}(k)$ thus we have
$R_{A}(k) = R_{B}(k)$ for $k
\le 2^{l} - 1$. This equality holds for every $l$, therefore we have
\begin{equation}
R_{A}(k) = R_{B}(k)\qquad \mbox{for every $k$}.
\end{equation}

To prove Theorem 2. and 3. we need the following three claims.

{\bf Claim 1.} Let
$0 < r_{1} < \dots{} < r_{s} \le m$ be integers. Then there exists at most
one pair of sets $(C, D)$ such that $C
\cup D = [0, m]$, $0 \in C$, $C \cap D = \{r_{1}, \dots{}, r_{s}\}$,
$R_{C}(k) = R_{D}(k)$ for every $k \le m$.

{\bf Proof of Claim 1.} We prove by contradiction. Assume that there exist at
least two pairs of different sets $(C_{1}, D_{1})$ and $(C_{2}, D_{2})$ which
satisfies the conditions of Claim 1. Let $v$ denote the smallest
positive integer such that $\chi_{C_{1}}(v) \ne \chi_{C_{2}}(v)$. It is
clear that $R_{C_{1}}(v) = R_{D_{1}}(v)$ and $R_{C_{2}}(v) =
R_{D_{2}}(v)$. We will prove that $R_{D_{1}}(v) = R_{D_{2}}(v)$ but
$R_{C_{1}}(v) \ne R_{C_{2}}(v)$ which is a contradiction. Obviously
\[
R_{D_{1}}(v) = |\{(d,d^{'}): d < d^{'}, d, d^{'} \in D_{1}, d + d^{'} = v\}|.
\]
As $0 \notin C \cap D$ and $0\in C$, we have $d$, $d^{'} < v$. We prove that $D_{1}
\cap [0, v - 1] = D_{2} \cap [0, v - 1]$, which implies that
$R_{D_{1}}(v) = R_{D_{2}}(v)$. Clearly we have $C_{1}
 \cap [0, v - 1] = C_{2} \cap [0, v - 1]$ and $[0, v - 1] = (C_{1}
 \cap [0, v - 1]) \cup (D_{1} \cap [0, v - 1])$ and $[0, v - 1] = (C_{2}
\cap [0, v - 1]) \cup (D_{2} \cap [0, v - 1])$. Let $(C_{1} \cap [0, v - 1]) \cap (D_{1} \cap [0, v
- 1]) = \{r_{1}, \dots{}, r_{t}\}$. Thus we have
$D_{1} \cap [0, v - 1] = \Big([0, v - 1] \setminus (C_{1} \cap [0, v -
1])\Big) \cup \{r_{1}, \dots{}, r_{t}\}$. Similarly $D_{2} \cap [0, v - 1] = \Big([0, v - 1] \setminus (C_{2} \cap [0, v -
1])\Big) \cup \{r_{1}, \dots{}, r_{t}\}$, which implies $D_{1} \cap [0, v - 1] = D_{2} \cap [0, v - 1]$. On the other hand as $0 \in C_{1} \cap C_{2}$ we have
\[
R_{C_{1}}(v) = |\{(c,c^{'}): c<c^{'} < v, c, c^{'} \in C_{1}, c + c^{'} =
v\}| + \chi_{C_{1}}(v), 
\]
and
\[
R_{C_{2}}(v) = |\{(c,c^{'}): c<c^{'} < v, c, c^{'} \in C_{1}, c + c^{'} =
v\}| + \chi_{C_{2}}(v),
\]
thus $R_{C_{1}}(v) \ne R_{C_{2}}(v)$.

{\bf Claim 2.} Let $(C,
D)$ be a pair of different sets, $C \cup D = [0, m]$, $C \cap D = \{r_{1},
\dots{}, r_{s}\}$, and $R_{C}(n) = R_{D}(n)$ for every $n$ nonnegative
integer and if $C^{'} = m - C$ and $D^{'} = m - D$ then $C^{'} \cup
D^{'} = [0, m]$, $C^{'} \cap D^{'} = \{m - r_{s},
\dots{}, m - r_{1}\}$, and $R_{C^{'}}(n) = R_{D^{'}}(n)$ for every positive 
integer $n$. 

{\bf Proof of Claim 2.} Clearly,
\[
R_{C}(k) = |\{(c,c^{'}): c < c^{'}, c, c^{'} \in C^{'}, c + c^{'} =
k\}| = |\{(c,c^{'}): c < c^{'}, m - c, m - c^{'} \in C, c + c^{'} =
k\}|
\]
\[
= |\{(m - c, m - c^{'}): c< c^{'}, m - c, m - c^{'} \in C, 2m - (c
+ c^{'}) = 2m - k\}| = R_{C}(2m - k).
\]
Similarly, $R_{D'}(k)=R_{D}(2m-k)$, which implies $R_{D'}(k)=R_{D}(2m-k)=R_{C}(2m-k)=R_{C'}(k)$, as desired.

{\bf Claim 3.} If for some positive integer $M$, the integers 
$M - 1, M-2,M-4,M-8,
\dots{}, M - 2^{u}$, $u=\lfloor log _2M\rfloor -1$ are all contained in the set 
$A$, then $M = 2^{u+1} - 1$.

{\bf Proof of Claim 3.} Suppose that the integers $M - 1, M-2,M-4,M-8,
\dots{}, M - 2^{u}$, $u=\lfloor log _2M\rfloor -1$ are all contained in the set $A$. If $M$ is even then $M - 2$ is also an even and $M-1=(M-2)+1$, therefore by the definition of $A$, 
$\chi _A(M-1)\ne \chi_A(M-2)$, thus we may assume
that $M$ is an odd positive integer. Assume that $M$ is not of the form $2^k-1$.
Obviously, $\chi_{A}(M) \ne \chi_{A}(M - 1)$.
Let $M = \sum_{i=0}^{u}b_{i}2^{i}$ be
the representation of $M$ in the number system with base $2$. Let $x$ denote
the largest index $i$ such that $b_{i} = 0$. Then $x\le \lfloor \log _2M\rfloor -1$. Thus we have
\[
M = \sum_{i=0}^{x-1}b_{i}2^{i} + 2^{x+1} + \sum_{i=x+2}^{u+1}b_{i}2^{i},\quad b_i\in \{ 0,1\},
\]
thus we have
\[
M - 2^{x} = \sum_{i=0}^{x-1}b_{i}2^{i} + 2^{x} +
\sum_{i=x+2}^{u+1}b_{i}2^{i} = \sum_{i=0}^{u+1}b^{'}_{i}2^{i},\quad b_i'\in \{ 0,1\},
\]
which implies that $\sum_{i=0}^{u}b_{i} = \sum_{i=0}^{u}b^{'}_{i}$. It follows from 
the definition of $A$ that
$\chi_{A}(M) = \chi_{A}(M - 2^{x})$. On the other hand
\[
\chi_{A}(M) \ne \chi_{A}(M - 1) = \chi_{A}(M - 2) = \chi_{A}(M - 4) = \dots{} = 
\chi_{A}(M - 2^{x}),
\]
which proves Claim 3.

{\bf Claim 4.} If for some positive integer $M$, the integers 
$M - 1, M-2,M-4,M-8,
\dots{}, M - 2^{u}$, $u=\lfloor log _2M\rfloor -1$ are all contained in the set 
$B$, then $M = 2^{u+1} - 1$. 

The proof of Claim 4. is completely the same as Claim 3.

Theorem 2. is a consequence of (1) and Claim 1. (for $s=0$).
In the next step we prove that if the sets $C$ and $D$ satisfies 
$C\cup D=[0,m]$, $C\cap D=\emptyset $ and $R_{C}(n) = R_{D}(n)$ for every 
positive integer 
$n$, then there exists an 
$l$ positive integer such that $C = A_{l}$ and $D = B_{l}$.
Then Claim 1. and (1) imply that $C = A
\cap [0, m]$ and $D = B \cap [0, m]$. Let $C^{'} = m - C$ and $D^{'} = m -
D$.  By Claim 1., Claim 2. and (1) we have $C^{'} = A
\cap [0, m]$ or $D^{'} = B \cap [0, m]$. It follows that
\[
\chi_{C^{'}}(2^{0}) = \chi_{C^{'}}(2^{1}) = \chi_{C^{'}}(2^{2}) = \dots{} = \chi_{C^{'}}(2^{u}),
\]
where $u= \lfloor \log _2M\rfloor -1$, which implies that
\[
\chi_{C}(m-1) = \chi_{C}(m-2) = \chi_{C}(m-4) = \dots{} = \chi_{C}(m-2^{u}).
\]
By Claim 3. and Claim 4. we get $m=2^{u+1}-1$. The proof of Theorem 3. is completed.

\section{Proof of Theorem 6.}

First, assume that there exists a positive integer $l$ such that $C = A_{2l}
\cup (2^{2l} - 1 + B_{2l})$, $D = B_{2l} \cup (2^{2l} - 1 +
A_{2l})$. Obviously, $C \cup D = [0, 2^{2l+1} - 2]$, $C \cap D = \{2^{2l}
- 1\}$, $0 \in C$ and we will prove that for every positive integer $n$, $R_{C}(n) = R_{D}(n)$. It is easy to see that

\[
R_{C}(n) = |\{(c,c^{'}): c < c^{'}, c, c^{'} \in A_{2l}, c + c^{'} =
n\}| +
\]
\[
|\{(c,c^{'}): c \in A_{2l}, c' \in 2^{2l} - 1 + B_{2l}, c + c' =
n\}|
\]
\[
+ |\{(c,c'): c < c', c, c' \in 2^{2l} - 1 + B_{2l}, c + c' =
n\}|
\]
\[
=  R_{A_{2l}}(n) + |\{(c,c'): c \in A_{2l}, c' \in B_{2l}, c + c' =
n - (2^{2l} - 1)\}| + R_{B_{2l}}(n - 2(2^{2l} - 1)).
\]
Moreover,
\[
R_{D}(n) = |\{(d,d'): d < d', d, d' \in B_{2l}, d + d' =
n\}| +
\]
\[
|\{(d,d'): d \in B_{2l}, d' \in 2^{2l} - 1 + A_{2l}, d + d' =
n\}|
\]
\[
+ |\{(d,d'): d, d' \in 2^{2l} - 1 + A_{2l}, d + d' =
n\}|
\]
\[
=  R_{B_{2l}}(n) + |\{(d,d'): d \in A_{2l}, d' \in B_{2l}, d + d' =
n - (2^{2l} - 1)\}| + R_{A_{2l}}(n - 2(2^{2l} - 1)).
\]

It follows from Theorem 3. that $R_{A_{2l}}(m)
= R_{B_{2l}}(m)$ for every positive integer $m$, thus we get the result.

In the next step we prove that if $C \cup D = [0, m]$, $C \cap D =
\{r\}$, $0 \in C$ and $R_{C}(n) = R_{D}(n)$, then there exists a positive integer $l$ such that $C = A_{2l}
\cup (2^{2l} - 1 + B_{2l})$, $D = B_{2l} \cup (2^{2l} - 1 + A_{2l})$ and
$m = 2^{2l + 1} - 2$, $r = 2^{2l} - 1$.
By Claim 2. if $C \cup D = [0, m]$, $C \cap D = \{r\}$ and $R_{C}(n) =
R_{D}(n)$ for every positive integer $n$, then for $C' = m - C$, $D^{'} = m - D$
 we have $C' \cup D^{'} = [0, m]$, $C' \cap D^{'} = \{m - r\}$ and by Claim 2. 
$R_{C'}(n) = R_{D^{'}}(n)$, thus we may assume that $r \le m/2$. 
Let
\begin{equation}
p_{C}(x) = \sum_{i=0}^{m}\chi_{C}(i)x^{i}, \hspace*{1mm} p_{D}(x) = \sum_{i=0}^{m}\chi_{D}(i)x^{i}.
\end{equation}
As $C \cup D = [0, m]$, $C \cap D = \{r\}$ we have 
\begin{equation}
p_{D}(x) = \frac{1 - x^{m+1}}{1 - x} - p_{C}(x) + x^{r}.
\end{equation}
Since $R_{C}(n) =
R_{D}(n)$ for every positive integer $n$, thus we have
\[
\sum_{n=0}^{\infty}R_{C}(n)x^{n} = \sum_{n=0}^{\infty}R_{D}(n)x^{n}.
\]
An easy observation shows that
\begin{equation}
\frac{1}{2}p_{C}(x)^{2} - \frac{1}{2}p_{C}(x^{2}) = \sum_{n=0}^{\infty}R_{C}(n)x^{n} = 
\sum_{n=0}^{\infty}R_{D}(n)x^{n} = \frac{1}{2}p_{D}(x)^{2} - \frac{1}{2}p_{D}(x^{2}).
\end{equation}
It follows from (3) and (4) that
\[
(p_{C}(x))^{2} - p_{C}(x^{2}) = \Big(\frac{1 - x^{m+1}}{1 - x} -
p_{C}(x) + x^{r}\Big)^{2} - \Big(\frac{1 - x^{2m+2}}{1 - x^{2}} -
p_{C}(x^{2}) + x^{2r}\Big).
\]
An easy calculation shows that
\begin{equation}
2p_{C}(x^{2}) = \frac{1 - x^{2m+2}}{1 - x^{2}} + 2p_{C}(x)\frac{1 -
x^{m+1}}{1 - x} - \Big(\frac{1 - x^{m+1}}{1 - x}\Big)^{2} -
2x^{r}\frac{1 - x^{m+1}}{1 - x} + 2x^{r}p_{C}(x).
\end{equation}

We will prove that $r$ must be odd. If $r$ were even and $r \le k \le 2r \le m$ is also even then it is easy to see from (2) and the coefficient of $x^k$ in (5) we have
\begin{equation}
2\chi_{C}\Big(\frac{k}{2}\Big) = 1 + 2\sum_{i \le k}\chi_{C}(i) - (k + 1) - 2
+ 2\chi_{C}(k - r),
\end{equation}
If $k + 1 \le 2r \le m$, then from the coefficient of $x^{k+1}$ in (5) we have
\begin{equation}
0 =  2\sum_{i \le k+1}\chi_{C}(i) -(k + 2) - 2
+ 2\chi_{C}(k + 1 - r).
\end{equation}
By (6) - (7) and dividing by 2 we get that
\[
\chi_{C}\Big(\frac{k}{2}\Big) = 1 - \chi_{C}(k + 1) + \chi_{C}(k - r)
- \chi_{C}(k + 1 - r).
\]
In view of Claim 1., $k+1-r<r$, $k-r$ is even, $C\cap [0,r-1]=A\cap [0,r-1]$ and by the definition of $A$ we get $\chi_{C}(k + 1 - r) + \chi_{C}(k - r) = 1$ thus we have
$\chi_{C}(k - r) - \chi_{C}(k + 1 - r) = \pm 1$. If $\chi_{C}(k - r)
- \chi_{C}(k + 1 - r) = 1$ then we get that $\chi_{C}(k - r) = 1$, $\chi_{C}(k
+ 1 - r) = 0$, which yields
 $\chi_{C}(k + 1) = 1$ and $\chi_{C}\Big(\frac{k}{2}\Big) = 1$.
On the other hand if $\chi_{C}(k - r)
- \chi_{C}(k + 1 - r) = -1$ then we get that $\chi_{C}(k - r) = 0$, $\chi_{C}(k
+ 1 - r) = 1$, which yields
 $\chi_{C}(k + 1) = 0$ and $\chi_{C}\Big(\frac{k}{2}\Big) = 0$.
This gives that $\chi_{C}\Big(\frac{k}{2}\Big) = \chi_{C}(k - r)$. Putting
$k = 2r - 2^{i+1}$, where $i+1\le \lfloor \log _2r\rfloor $ we obtain that
$\chi_{C}(r - 2^{i+1}) = \chi_{C}(r - 2^{i})$. It follows that
\[
\chi_{C}(r - 1) = \chi_{C}(r - 2) = \chi_{C}(r - 4) = \dots{} = \chi_{C}(r - 2^{t}).
\]
for $t=\lfloor \log _2r\rfloor -1$. Moreover, Claim 3. and Claim 4. imply that $r=2^l-1$, which contradicts the assumption that $r$ is even.

It follows that $r$ must be odd. By using the same argument as before we prove that $r = 2^{2l} - 1$. If $r \le k < 2r \le m$ and $k$ is even then from the coefficient of $x^k$ in (5) we have
\[
2\chi_{C}\Big(\frac{k}{2}\Big) = 1 + 2\sum_{i \le k}\chi_{C}(i) - (k + 1) -
2 + 2\chi_{C}(k-r).
\]
In this case $k - 1$ is odd, and $k - 1 \ge r$, therefore from the coefficient of $x^{k-1}$ in (5) we have
\[
0 = 2\sum_{i \le k-1}\chi_{C}(i)-k  - 2 + 2\chi_{C}(k-1-r).
\]
As before, subtracting the above equalities and dividing by 2 we get that
\[
\chi _{C}\Big(\frac{k}{2}\Big) = \chi_{C}(k) + \chi_{C}(k-r) - \chi_{C}(k-1-r).
\]
If $r$ is odd, then $k-1-r$ is even, we know from Claim 1. that $C\cap [0,r-1]=A\cap [0,r-1]$ and by definition of $A$ we get $\chi_{C}(k-1-r) + \chi_{C}(k-r) = 1$ thus we have
$\chi_{C}(k-r) - \chi_{C}(k-1-r) = \pm 1$. If $\chi_{C}(k-r) - \chi_{C}(k-1-r) =
1$ then we get that $\chi_{C}(k-r) = 1$, $\chi_{C}(k-1-r) = 0$, which yields
 $\chi_{C}(k) = 0$ and $\chi_{C}\Big(\frac{k}{2}\Big) = 1$.

If $\chi_{C}(k-r) - \chi_{C}(k-1-r) = -1$ then we get that $\chi_{C}(k-r) = 0$, $\chi_{C}(k-1-r) = 1$, which yields
 $\chi_{C}(k) = 1$ and $\chi_{C}\Big(\frac{k}{2}\Big) = 0$. This gives
 that $\chi_{C}(k-r) = \chi_{C}\Big(\frac{k}{2}\Big)$ when $r \le k <
 2r$ and $k$
 is even. Put $k = 2r - 2^{i+1}$, where $i+1\le \lfloor \log _2r\rfloor $ implies that $\chi_{C}(r - 2^{i+1}) = \chi_{C}(r
 - 2^{i})$. It follows that
\[
\chi_{C}(r - 1) = \chi_{C}(r - 2) = \chi_{C}(r - 4) = \dots{} = \chi_{C}(r - 2^{\lfloor \log _2r\rfloor -1}),
\]
which yields by Claim 3. and Claim 4. that $r = 2^{u} - 1$. If $k = r$,
then from the coefficient of $x^{k}$ in (5) we have
\begin{equation}
0 =  2\sum_{i \le r}\chi_{C}(i) - (r + 1) - 2 + 2.
\end{equation}
On the other hand if  $k = r - 1$, then from the coefficient of $x^{k-1}$ in (5) we have
\begin{equation}
2\chi_{C}\Big(\frac{r - 1}{2}\Big) =  1 + 2\sum_{i \le r-1}\chi_{C}(i) - r .
\end{equation}
Since $r \in C$, it follows from (8) - (9) that
\[
-2\chi_{C}\Big(\frac{r - 1}{2}\Big) =  -2 + 2\chi_{C}(r) = 0,
\]
we get that $0 = \chi_{C}\Big(\frac{r - 1}{2}\Big) = \chi_{A}\Big(\frac{r - 1}{2}\Big) $, thus $\frac{r - 1}{2} =
2^{u-1} - 1$, where $u - 1$ is odd, so that $r = 2^{2l} - 1$.

As $r\le m/2$, Claim 1. and the first part of Theorem 3. imply that $C \cap [0, 2r] = A_{2l} \cup (2^{2l} - 1 + B_{2l})$ and $D \cap
[0, 2r] = B_{2l} \cup (2^{2l} - 1 + A_{2l})$.

We will show that $m < 3 \cdot 2^{2l} - 2$. We prove by
contradiction. Assume that $m > 3 \cdot 2^{2l} - 2$. We will prove that
\[
C \cap [0,3\cdot 2^{2l}-3]= A_{2l} \cup (2^{2l} - 1 + B_{2l}) \cup (2^{2l+1} - 1 + (B_{2l} \cap [0,
2^{2l} - 2]))
\]
and
\[
D \cap [0,3\cdot 2^{2l}-3] = B_{2l} \cup (2^{2l} - 1 + A_{2l}) \cup (2^{2l+1} - 1 + (A_{2l} \cap [0,
2^{2l} - 2])).
\]
Define the sets $E$ and $F$ by 
\[
E = A_{2l} \cup (2^{2l} - 1 + B_{2l}) \cup (2^{2l+1} - 1 + (B_{2l} \cap [0, 2^{2l} - 2])),
\]
and
\[
F = B_{2l} \cup (2^{2l} - 1 + A_{2l}) \cup (2^{2l+1} - 1 + (A_{2l} \cap [0, 2^{2l} - 2])).
\]
In this case $R_{E}(m) = R_{F}(m)$ for every $m \le 2^{2l+1} - 2$.
If $2^{2l+1} -
2 < n \le  3 \cdot 2^{2l} - 3$, then we have
\[
R_{E}(n) = |\{(c,c'): c \in A_{2l}, c' \in 2^{2l} - 1 + B_{2l},  c +
c' = n\}| +
\]
\[
|\{(c,c'): c \in A_{2l}, c' \in 2^{2l+1} - 1 +
(B_{2l} \cap [0, 2^{2l} - 2]), c + c' = n\}| +
\]
\[
|\{ (c,c'): c<c',c,c'\in 2^{2l}-1+B_{2l}, c+c'=n\}|
\]
\[
= |\{(c,c'): c \in A_{2l}, c' \in B_{2l}, c + c' = n - (2^{2l} - 1)\}|
\]
\[
+ |\{(c,c'): c \in A_{2l}, c' \in B_{2l}, c + c' = n - (2^{2l+1} -
1)\}| +
\]
\[
R_{B_{2l}}(n-2(2^{2l}-1))
\]
and
\[
R_{F}(n) = |\{(d,d^{'}): d \in B_{2l}, d^{'} \in 2^{2l} - 1 + A_{2l}, d +
d' = n\}|
\]
\[
+ |\{(d,d'): d \in B_{2l}, d' \in 2^{2l+1} - 1 +
(A_{2l} \cap [0, 2^{2l} - 2]), d + d' = n\}|
\]
\[
||\{
(d,d'):d<d',d,d'\in 2^{2l}-1+A_{2l},d+d'=n
\}
\]
\[
=|\{(d,d'): d \in B_{2l}, d' \in A_{2l}, d + d' = n
- (2^{2l} - 1)\}|
+
\]
\[
|\{(d,d'): d \in B_{2l}, d' \in A_{2l}, d + d' = n - (2^{2l+1} -
1)\}|+R_{A_{2l}}(n-2(2^{2l}-1)),
\]
which imply $R_E(n)=R_F(n)$ for every positive integer $n$ by Theorem 3., therefore by Claim 1., this is the only possible starting of $C$ and $D$.

We will prove
that $3 \cdot 2^{2l} - 2 \in C$. We prove by contradiction. Assume that
$3 \cdot 2^{2l} - 2 \in D$, that is 
\[
C \cap [0,3\cdot 2^{2l}-2]=A_{2l} \cup (2^{2l} - 1 + B_{2l}) \cup (2^{2l+1} - 1 + B_{2l})
\]
and 
\[
D \cap [0,3\cdot 2^{2l}-2]=B_{2l} \cup (2^{2l} - 1 + A_{2l}) \cup (2^{2l+1} - 1 + A_{2l}). 
\]
We have a solution $3 \cdot 2^{2l} - 2 =
(2^{2l+1} - 1) + (2^{2l} - 1)$ in set $D$, thus we have
\[
R_{D}(3 \cdot 2^{2l} - 2) = 1 + |\{(d,d'): d \in B_{2l}, d' \in 2^{2l+1} - 1 +
A_{2l}, d + d' = 3 \cdot 2^{2l} - 2\}|
\]
\[
+ |\{ (d,d'):d<d', d,d'\in 2^{2l}-1+A_{2l},d+d'=3\cdot 2^{2l}-2 \}|
\]
\[
= 1 + |\{(d,d'): d \in B_{2l}, d' \in A_{2l}, d + d' = 2^{2l} - 1)\}|+R_{A_{2l}}(2^{2l}).
\]
On the other hand
\[
R_{C}(3 \cdot 2^{2l} - 2) = |\{(c,c'): c \in A_{2l}, c' \in 2^{2l+1}
- 1 + B_{2l}, c + c' = 3 \cdot 2^{2l} - 2\}|
\]
\[
+ |\{ (c,c'):c<c', c,c'\in 2^{2l}-1+B_{2l},c+c'=3\cdot 2^{2l}-2 \}|
\]
\[
=|\{(c,c'): c \in A_{2l}, c' \in B_{2l}, c + c' = 2^{2l} - 1)\}| + R_{B_{2l}}(2^{2l}),
\]
therefore by Theorem 3 we have $R_{D}(3 \cdot 2^{2l} - 2)>R_{C}(3 \cdot 2^{2l} - 2)$, which is a contradiction. We may assume that $3 \cdot 2^{2l} - 2 \in C$.

Using the fact $1\not \in C$ we have
\[
R_{C}(3 \cdot 2^{2l} - 1) = |\{(c,c'): c \in A_{2l}, c' \in 2^{2l+1}
- 1 +
B_{2l}, c + c' = 3 \cdot 2^{2l} - 1\}|
\]
\[
+|\{ (c,c'):c<c',c,c'\in 2^{2l}-1+B_{2l},c+c'=3\cdot 2^{2l}-1 \}| + \chi_{C}(3 \cdot 2^{2l} - 1)
\]
\[
= |\{(c,c'): c \in A_{2l}, c' \in B_{2l}, c + c' = 2^{2l} \}| + R_{B_{2l}}(2^{2l}+1)+ 
\chi_{C}(3 \cdot 2^{2l} - 1).
\]
On the other hand using $1\in B_{2l}$, $2^{2l}-1\in A_{2l}$ and $3\cdot 2^{2l}-2\in C$ we get
\[
R_{D}(3 \cdot 2^{2l} - 1) = |\{(d,d'): d \in B_{2l}, d' \in 2^{2l+1}
- 1 +
A_{2l}\cap[0,2^{2l}-2], d + d' = 3 \cdot 2^{2l} - 1\}| 
\]
\[
+ |\{ (d,d'):d<d',d,d'\in 2^{2l}-1+A_{2l},d+d'=3\cdot 2^{2l}-1 \}|
\]
\[
= |\{(d,d'): d \in B_{2l}, d' \in  A_{2l}, d + d' = 2^{2l} \}| +R_{A_{2l}}(2^{2l}+1)-1,
\]
 therefore by Theorem 3 we have $R_{C}(3 \cdot 2^{2l} - 1)>R_{D}(3 \cdot 2^{2l} - 1)$, which is a contradiction, that is we have $m < 3 \cdot 2^{2l} - 2$.

It follows that
\[
C = A_{2l} \cup (2^{2l} - 1 + B_{2l}) \cup (2^{2l+1} - 1 + (B_{2l} \cap [0, m - (2^{2l+1} - 1)]))
\]
and
\[
D = B_{2l} \cup (2^{2l} - 1 + A_{2l}) \cup (2^{2l+1} - 1 + (A_{2l} \cap [0, m - (2^{2l+1} - 1)])).
\]

We will prove that $m = 2^{2l+1} - 2$. Assume that $m > 2^{2l+1} - 2$. If $m - (2^{2l+1} -
1) \ne 2^{k} - 1$, then by Corollary 1., there exists an $m - (2^{2l+1} - 1) < u < 2(m -
(2^{2l+1} - 1))$ such that
\[
R_{A \cap [0, m - (2^{2l+1} - 1)]}(u) \ne R_{B \cap [0, m - (2^{2l+1} -
1)]}(u).
\]
Since $m+2^{2l+1}-1<u+2(2^{2l+1}-1)<2m$ we obtain
\[
R_{C}(2(2^{2l+1} - 1) + u) = R_{2^{2l+1} - 1 + B_{2l} \cap [0, m - (2^{2l+1} -
1)]}(2(2^{2l+1} - 1) + u) = R_{B_{2l} \cap [0, m - (2^{2l+1} - 1)]}(u).
\]
Similarly
\[
R_{D}(2(2^{2l+1} - 1) + u) = R_{A_{2l} \cap [0, m - (2^{2l+1} - 1)]}(u),
\]
which is a contradiction. Thus we may assume that $m - (2^{2l+1} - 1) = 2^{k} - 1$, where $k <
2l$. Hence
\[
C = A_{2l} \cup (2^{2l} - 1 + B_{2l}) \cup (2^{2l+1} - 1 + B_{k})
\]
and
\[
D = B_{2l} \cup (2^{2l} - 1 + A_{2l}) \cup (2^{2l+1} - 1 + A_{k}).
\]
If $k=0$, then $C,D\subset [0,2^{2l+1}-1]$ and $2^{2l+1}-2,2^{2l+1}-1\in D$, therefore $R_{C}(2^{2l+2}-3)=0$ and $R_{D}(2^{2l+2}-3)=1$, a contradiction. Thus we may assume that $k>0$. Then if $C' = 2^{2l+1} + 2^{k} - 2 - C$ and $D' = 2^{2l+1} + 2^{k}
- 2 - D$ it follows that $C' \cup D' = [0, 2^{2l+1} + 2^{k} - 2]$,
$C' \cap D' = \{2^{2l} + 2^k - 1\}$ and by Claim 2. $R_{C'}(n) = R_{D'}(n)$ for every positive integer $n$.
Thus we have $C' \cap [0, 2^{2l} + 2^{k} - 2] = A \cap [0, 2^{2l} +
2^{k} - 2]$ or $C' \cap [0, 2^{2l} + 2^{k} - 2] = B \cap [0, 2^{2l} +
2^{k} - 2]$. We prove that $C'\cap [0,2^{2l}+2^k-2]=A\cap [2^{2l}+2^k-2]$. Assume that $C'\cap [0,2^{2l}+2^k-2]=B\cap [2^{2l}+2^k-2]$. Then $C'\cap [0,2^k-1]=B_k$ and $C'\cap [2^k,2^{k+1}-1]=2^k+A_k$, which is a contradiction because $C'\cap [2^k,2^{k+1}-1]=2^k+B_k$. We know $C' \cap [0, 2^{k} - 1] = A_{k}$, then $C' \cap [2^{k},
2^{k+1} - 1] = 2^{k} + B_{k}$, and $C' \cap [2^{k+1}, 3 \cdot 2^{k} - 2] =
2^{k+1} + B_{k}\cap [0,2^k-2]$. On the other hand for $C^{''} = 2^{2l+1} - 2 - (A_{2l} \cup (2^{2l}
- 1 + B_{2l}))$ we have $C''\cap [0,2^{2l}-2]=B\cap [0,2^{2l}-2]$, therefore $C^{''}\cap [0,2^k-1]=B_k$ and $C^{''}\cap [2^k,2^{k+1}-2]=2^k+A_k\cap [0,2^k-2]$, which is a contradiction because $C'=A_k\cup (2^k+C^{''})$.
 The proof of Theorem 6. is completed.

\end{document}